\documentclass[reqno]{amsart}

\usepackage{amsmath}%
\usepackage{amsfonts}%
\usepackage{amssymb}%
\usepackage{ifpdf}
\ifpdf
\usepackage{graphicx}
\else
\usepackage[draft]{graphicx}
\fi
\usepackage{geometry}
\usepackage{subfig}
\newtheorem{theorem}{Theorem}
\theoremstyle{plain}

\newcommand {\Log}{\text{Log}\,}
\begin{document}

\title[Prime Counting Function]{The Analytic Expression for Riemann's Prime Counting Function via the Residue Theorem}
\author{Dominic C. Milioto} 

\subjclass[2000]{11M06,11-02 (Primary) 30-02,44A10 (Secondary)}
\keywords{Prime counting function, Residue Theorem, Inverse Laplace Transform, zeta function, argument function, holomorphic branch}

\begin{abstract}
In his paper "On the Number of Primes Less Than a Given Magnitude"\cite{Riemann}, Bernhard Riemann introduced a prime counting function $F(x)$ which counts the number of primes under $x$.  Riemann obtained an analytic expression for $F(x)$ by evaluating an inverse Laplace Transform.  His method involved advanced techniques of analysis.  However, this transform can be evaluated using the Residue Theorem when an appropriate branch of $\log\zeta(s)$ is defined.  In this paper, a method for constructing a holomorphic branch of $\log\zeta(s)$ extending to the left half-plane is described along with it's geometry surrounding the logarithmic branch points.  Using this information, an integral representation of $F(x)$ is formulated in terms of this branch of $\log\zeta$ which is then evaluated.  The results are shown equal to Riemann's expression.
\end{abstract}

\maketitle

\section{Introduction}

The focus of this paper is three-fold:  
\begin{enumerate}
\item To review an historically significant problem in Analytic Number Theory,
\item To describe a particular branch of $\log\zeta(s)$,
\item Using (2), solve (1) in a more elementary way.
\end{enumerate}

In 1859, Bernhard Riemann published ``On the Number of Primes Less than a Given Magnitude''$\footnote{See Edwards for an English translation of Riemann's paper.}$, and introduced an analytic continuation of the Euler product he called the zeta function: $\zeta(s)$.  Using $\zeta(s)$, infinite product representations of the $\Gamma$ and $\xi$ functions and relying on advanced principles of analysis, Riemann went on to derive an analytic expression for a function he called $F(x)$ equal to the number of prime numbers under $x$.  However, if a branch of $\log\zeta(s)$ is carefully chosen, Riemann's expression for $F(x)$ can be derived via the Residue Theorem using simpler techniques of analysis.  This approach is described below.

Throughout this paper, a complex variable is represented as $s=\sigma+it$, and the non-trivial zeros of the zeta function with positive imaginary component as $\rho_n=\sigma_n+i\gamma_n$ (or just $\rho$ in some cases), such that $\gamma_n<\gamma_{n+1}$.  The principal branch of the complex logarithm is denoted by Log and other branches are given by log.  Common in Complex Analysis, Log represents logarithm to base e, that is Log(e)=1.  That convention is used here.  In all cases, the variable x is taken to be a real number greater than one. 

Riemann derives his prime counting function by beginning with the Euler product:

$$\zeta (s)=\prod _p \frac{1}{1-p^{-s}},\quad \text{Re}(s)>1.$$
Taking Logarithms of the above gives
\begin{equation}
\text{Log}\;\zeta (s)=-\sum _p \text{Log}\left(1-p^{-s}\right),
\label{firsteqn}
\end{equation}
and remembering $\text{Log}(1-z)= -\sum _{\text{\textit{$n$}}=1}^{\infty } \frac{\text{\textit{$z$}}^{\text{\textit{$n$}}}}{\text{\textit{$n$}}}$ for $|$z$|<$1, we obtain

$$\text{Log}\;\zeta (s)=\sum_{p}\sum_{n=1}^{\infty} \frac{1}{n} p^{-ns},\quad \text{Re}(s)>1.$$
Replacing each prime power with
$$p^{-s}=s\int_p^{\infty } x^{-s-1} \, dx,\quad p^{-2s}=s\int_{p^2}^{\infty } x^{-s-1} \, dx,\cdots,p^{-ns}=s\int_{p^n}^{\infty} x^{-s-1}dx,\cdots,$$
we can express Equation \eqref{firsteqn} as

\begin{equation}
\frac{\text{Log}\; \zeta (s)}{s}=\int _1^{\infty }f(x)x^{-s-1}dx
\label{pcf}
\end{equation}
with
\begin{equation*}
f(x)=F(x)+\frac{1}{2} F(x^{1/2})+\frac{1}{3} F(x^{1/3})+\cdots+\frac{1}{n}F(x^{1/n})+\cdots,\quad n=1,2,3,\cdots.
\end{equation*}
$F(x)$ is equal to the number of primes strictly less than $x$ and as is common in Fourier Theory, at points of discontinuity, namely at prime numbers, it is defined to be one-half the difference between it's old and new value.  With this definition of $F(x)$, we see $f(x)$ begins at 0 for x=0 and increases by a jump of 1 at primes, by a jump of 1/2 at prime squares $p^2$, by a jump of $1/3$ at prime cubes $p^3$, etc.  That is, $f(x)$ is zero for $0\leq x<2$, is $1/2$ for $x=2$, is 1 for $2<x<3$, is $1\frac{1}{2}$ for $x=3$, is $2$ for $3<x<4$, is $2\frac{1}{4}$ for $x=4$, is $2\frac{1}{2}$ for $4<x<5$, is $3$ for $x=5$ and so on.
From $f(x)$, one can find $F(x)$ via M\"{o}bius inversion to obtain
\begin{equation*}
F(x)=f(x)-\frac{1}{2} f(x^{1/2})-\frac{1}{3}f(x^{1/3})+\frac{1}{6}f(x^{1/6})+\cdots+\frac{\mu(n)}{n} f(x^{1/n})
\end{equation*}
with $n$ a positive integer such that $n\leq\frac{\text{Log}(x)}{\text{Log}2}$.  The function F(x) is Riemann's prime counting function and Riemann's derivation of an analytic expression for $f(x)$ and thus $F(x)$, is the principal result of his paper.  He accomplished this by inverting \eqref{pcf} via Fourier inversion giving

\begin{equation}
f(x)=\frac{1}{2\pi  i}\underset{a-i \infty }{\overset{a+i \infty }{\int }}\frac{\text{Log}\; \zeta (s)}{s}x^s ds,\quad a>1, x>1.
\label{eqn002}
\end{equation}
In order to evaluate \eqref{eqn002}, Riemann expressed $\zeta(s)$ in terms of the $\xi $ function:

$$\xi (s)=1/2 s(s-1)\pi ^{-s/2}\Gamma (s/2)\zeta (s),\quad s\in \mathbb{C},$$
and solved for $\text{Log}\,\zeta$.  He then expressed both $\xi$ and $1/\Gamma$ in terms of their infinite product representations and proceeded to evaluate the resulting integrals  directly to obtained \footnote{Riemann's paper contains an error in the expression for $f(x)$.  $(\ref{eqn003})$ is the corrected version.}

\begin{equation}
f(x)=Li(x)-\sum _{Im(\rho )>0} \left[Li\left(x^{\rho }\right)+Li\left(x^{1-\rho }\right)\right]+\int_x^{\infty } \frac{1}{t\left(t^2-1\right) \text{Log}(t)} \ dt-\text{Log}\,(2)
\label{eqn003}
\end{equation}
in which $Li(x)$ is the logarithmic integral and 
\begin{equation}
Li(x^{\rho})=P.V.\int_0^x \frac{s^{\rho-1}}{\Log(s)}ds=
\begin{cases}\displaystyle
\mathop\int\limits_{-\infty+iv}^{u+iv}\frac{e^s}{s}ds+\pi i, & \text{Im}(\rho)>0  \vspace{10pt} \\
\displaystyle
\mathop\int\limits_{-\infty+iv}^{u+iv}\frac{e^s}{s}ds-\pi i, & \text{Im}(\rho)<0,  \\
\end{cases}
\label{ingham001}
\end{equation}
with $u+iv=\rho\ln(x).$  The exponential expressions are obtained by first considering
$$
\mathop\int\limits_{C}\frac{s^{\rho-1}}{\Log(s)}ds=P.V.\int_{0}^{x}\frac{s^{\rho-1}}{\Log(s)}ds\pm \pi i
$$
where the contour $C$ is along the real axis from zero to $x$ with an indentation above the (simple) pole at $s=1$ for $\text{Im}(\rho)>0$ giving rise to the $+\pi i$ term,  or below it for $\text{Im}(\rho)<0$ and giving rise to the $-\pi i$ term.  Letting $s=e^{u}$ gives
$$
\mathop\int\limits_{C}\frac{s^{\rho-1}}{\Log(s)}ds=\int_{-\infty}^{\Log(x)}\frac{e^{u\rho}}{u}du
$$ 
where the path of integration for the integral on the right is now along the real axis with an indentation above or below the origin.  Making the substitution $v=u\rho$, then
$$
\mathop\int\limits_{C}\frac{s^{\rho-1}}{\Log(s)}ds=\int_{-\infty}^{\Log(x)}\frac{e^{u\rho}}{u}du=\mathop\int\limits_{-\infty-i\infty}^{\rho\ln(x)} \frac{e^{v}}{v}dv
$$
in which
$$\mathop\int\limits_{-\infty-i\infty}^{\rho\ln(x)} \frac{e^{v}}{v}dv=\mathop\int\limits_{-\infty-i\infty}^{-\infty+i\gamma_n\ln(x)} \frac{e^{v}}{v}dv+\mathop\int\limits_{-\infty+i\gamma_n\ln(x)}^{\rho\ln(x)} \frac{e^{v}}{v}dv.
$$
The center integral of the last expression tends to zero leaving the desired result.  For a further explanation of \eqref{ingham001}, see Edwards and Ingham \cite{Edwards,Ingham}.

Although Riemann's work represents a beautiful case study in analysis, his approach to evaluating \eqref{eqn002} was based on advanced methods of analysis.  In this paper,  \eqref{eqn002} is evaluated more simply using a holomorphic branch of $\log\zeta(s)$.  The geometry of this branch is investigated and the results used to derive an integral expression for $f(x)$ which is then  evaluated via the Residue Theorem.

\section{Constructing the holomorphic extension of $\text{Log}\,\zeta(s)$ to the left half-plane}
\label{section2}
\begin{figure}
	\centering
		\includegraphics{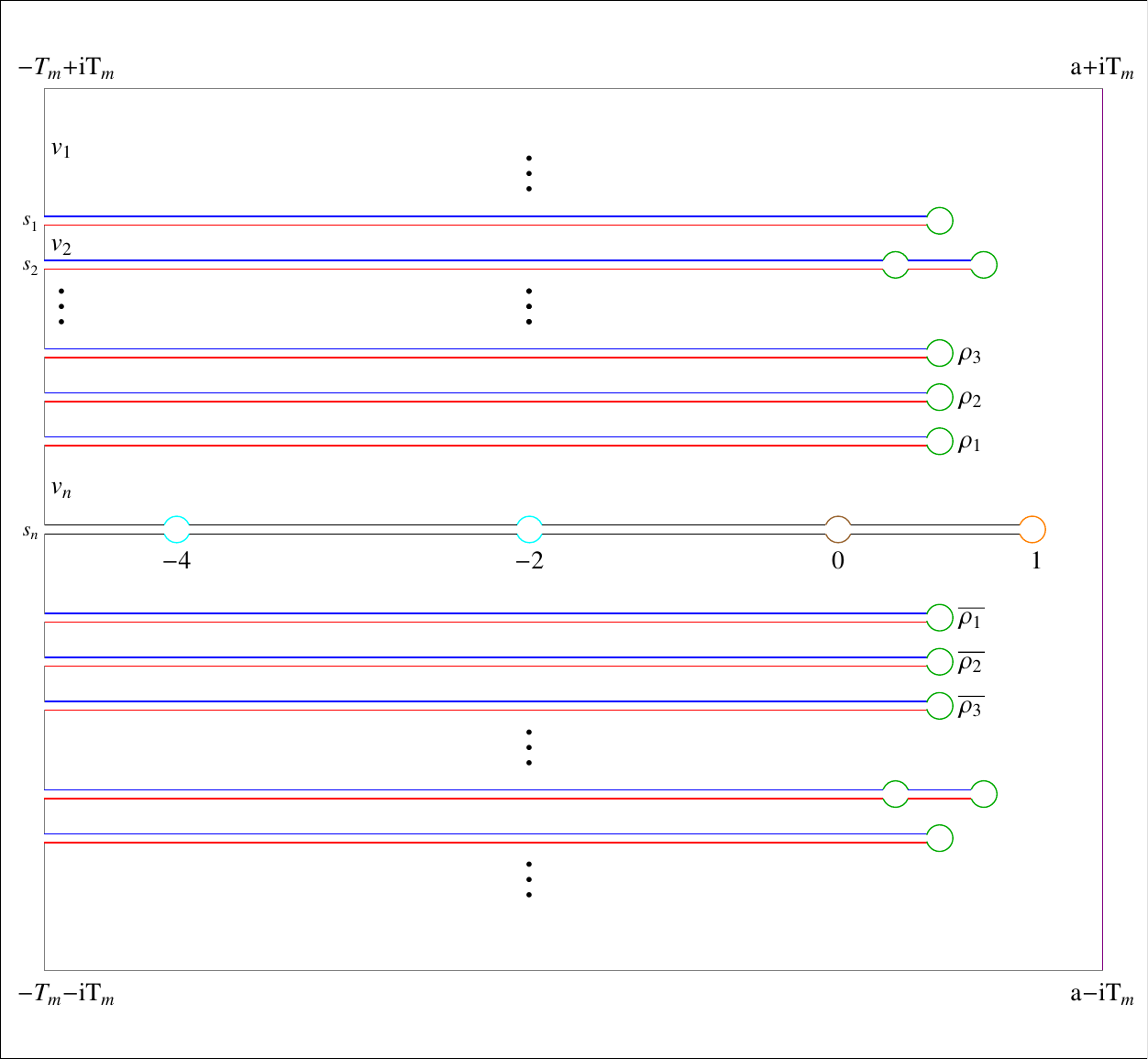}
		\caption{B(m) contour}
   \label{fig:file1}
\end{figure}

The function $\text{Log}\,\zeta(s)$ is holomorphic in the half-plane $\text{Re}(s)>1$ since $\zeta(s)$ is non-zero there.  However, there exists a holomorphic branch that extends to the left half-plane.  Proof of the following theorem can be found in \cite{Stein}.
\begin{theorem}
If $f$ is a nowhere vanishing holomorphic function in a simply connected region $\Omega$, then there exists a holomorphic function $g$ on $\Omega$ such that
$$f(z)=e^{g(z)}.$$
The function $g(z)$ is given by $\log f(z)$, and determines a ``branch'' of that logarithm.  Also:
$$g'(z)=\frac{f'(z)}{f(z)}.$$
\label{theorem001}
\end{theorem}
The trivial zeros of $\zeta(s)$ are located along the negative real axis at the points $s=-2n$.  The non-trivial zeros are symmetric in the critical strip $0<\sigma<1$ to the line $\sigma=1/2$, and $\zeta$ has one simple pole at $s=1$.  Therefore, we can define $\Omega$ to be the complex plane with branch cuts at the real axis extending from $-\infty$ to one, and branch cuts of the form $\sigma\pm i\gamma_n$ for each zero on the critical line and branch cuts between zeros for those off the critical line.  We then construct a contour $B(m)$ around $\Omega$.  This contour is shown in Figure (\ref{fig:file1}) and explained below.  Then $\log\zeta(s)=\ln|\zeta(s)|+i\arg \zeta(s)$ is holomorphic in $\Omega$ with

\begin{equation*}
\frac{d}{ds}\log\zeta(s)=\frac{\zeta'(s)}{\zeta(s)}, \quad s\notin \mathbb{B},
\end{equation*}
where the set $\{\mathbb{B}\}$ represents the branch cuts, and the function $\arg$ is now a continuous and analytic function of $s$ throughout $\Omega=\mathbb{C}\backslash\{\mathbb{B}\}$. 
Therefore,
\begin{equation}
\log\zeta(s)=\int_{s_0}^s\frac{\zeta'(s)}{\zeta(s)}ds+\log\zeta(s_0)
\label{inteqn1}
\end{equation}
with the path of integration remaining in $\Omega$.

%
% Figure of contour
%

In order use $\log\zeta(s)$ in this analysis, we must first understand the geometry of $\arg\zeta(s)$ which we do now.

\section{Local Geometry of $\log\zeta(s)$} 
\label{sectiona}

$\text{Log}\,|\zeta(s)|$ is continuous except at the zeros of the zeta function and the point $s=1$.  The imaginary sheet, $\arg\zeta(s)$, is however more complicated, and a description of this sheet in the neighborhood of the branch cuts is needed in order to use the Residue Theorem to derive Riemann's expression for $f(x)$. 

\subsection{Geometry of $\log\zeta(s)$ along real axis}
\label{section31}

\begin{figure}
	\centering
		\includegraphics{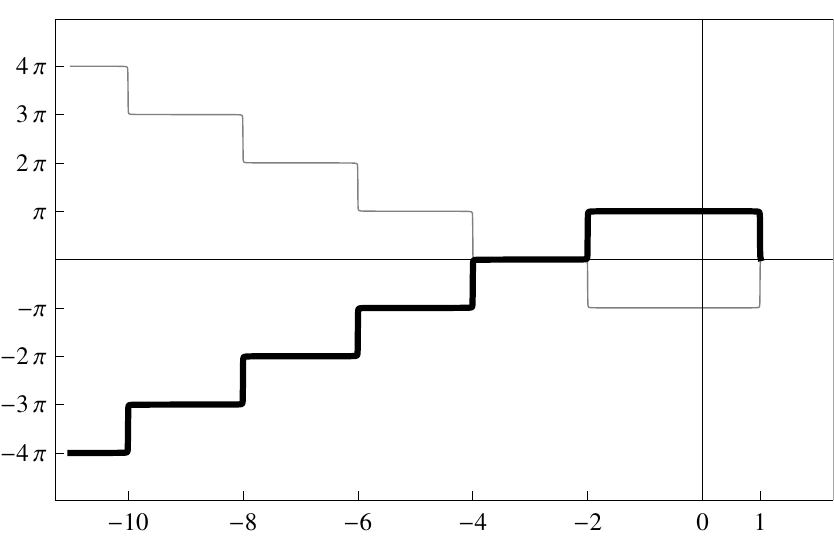}
	\caption{$\arg\zeta(s)$ along black contours of B(m)}
	\label{fig:realaxiscloseup}
\end{figure}

\begin{figure}
	\centering
		\includegraphics{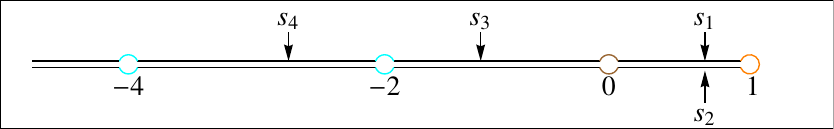}
	\caption{Close up of black contour at $s=-2$}
	\label{fig:closeupat2}
\end{figure}

$\arg\zeta(s)$ along a section of the black contours of $B(m)$ is shown in Figure (\ref{fig:realaxiscloseup}). Note the stair-step geometry.  The dark black trace is the contour below the real axis and the gray trace, above the real axis .  Empirical results suggest the difference in argument changes by $2\pi$ at each step. We can prove this difference using the following Theorem \cite{Marsden}:

\begin{theorem}
Let $f(s)$ be analytic with a simple pole at $s_0$ and $\Gamma_r$ be an arc of a circle of radius $r$ and angle $\alpha$ centered at $s_0$.  Then:
\begin{equation}
\lim_{r\to 0}\int_{\Gamma_r} f(s)ds=\alpha i\mathop\textnormal{Res}_{s=s_0} f(s).
\end{equation}
\label{theorem005}
\end{theorem}
Consider now
\begin{equation*}
\log\zeta(s)=\int_{s_0}^s\frac{\zeta'(s)}{\zeta(s)}ds+\log\zeta(s_0)
\end{equation*}
and the points $s_1$ and $s_2$ in Figure (\ref{fig:closeupat2}) where $s_1=\sigma+i r$ and $s_2=\sigma-i r$ with $-2<\sigma<1$ (note the indentation around the origin applies only to $\frac{\log\zeta(s)}{s}$).  We then have in the limit as the black contours approach the real axis,

\begin{equation}
\log\zeta(s_2)\to \lim_{r\to 0}\left\{\int_{\sigma+ir}^{1-r+i r}\frac{\zeta'(s)}{\zeta(s)}ds+\mathop\int\limits_{\Gamma_r}\frac{\zeta'(s)}{\zeta(s)}ds-\int_{\sigma-ir}^{1-r-ir}\frac{\zeta'(s)}{\zeta(s)}ds\right\}+\log\zeta(s_1).
\label{eqn33}
\end{equation}
Now, the first and third integrals are analytic outside some deleted neighborhood $D(1,\delta)$ of $s=1$ and therefore the integrals in this region cancel as $r\to 0$.  Inside $D(1,\delta)$,

\begin{equation}
\frac{\zeta'(s)}{\zeta(s)}=\frac{-1}{s-1}+\frac{\phi'(s)}{\phi(s)}
\label{eqn805}
\end{equation}
with $\phi(s)\ne 0$ for $s\in D(1,\delta)$.  Clearly, the residue of $\frac{\zeta'(s)}{\zeta(s)}$ at $s=1$ is $-1$. 

Now, because of the symmetric form of the limit, i.e., the terms $1-r\pm ir$, the angle $\alpha$ of $\Gamma_r$ will always extend from $3\pi/4$ to $-3\pi/4$ even in the limit as $r\to 0$.  Therefore
$$\lim_{r\to 0}\mathop\int\limits_{\Gamma_r}\frac{\zeta'(s)}{\zeta(s)}ds=\frac{3\pi}{2} i.$$

Consider now the remaining two integrals inside $D(1,\delta)$.  Substituting  \eqref{eqn805} into the integrals, the $\frac{\phi'(s)}{\phi(s)}$ terms will again because they are analytic, cancel as $r\to 0$ leaving
\begin{equation*}
\begin{aligned}
&\lim_{r\to 0}\left\{\int_{\sigma+ir}^{1-r+ir}\frac{-1}{s-1}ds+\int_{1-r-ir}^{\sigma-ir}\frac{-1}{s-1}ds\right\} \\
&=-\lim_{r\to 0}\left\{\int_{\sigma+ir}^{1-r+ir}\frac{1}{s-1}ds-\int_{\sigma-ir}^{1-r-ir}\frac{1}{s-1}ds\right\} \\
&=-\lim_{r\to 0}\Big\{\Log(-r+ir)-\Log(\sigma-1+ir)-\Log(-r-ir)+\Log(\sigma-1-ir)\Big\} \\
&=-\Big(\frac{3\pi}{4} i-\pi i+\frac{3\pi}{4} i-\pi i\Big)=\frac{\pi}{2}i.
\end{aligned}
\end{equation*}

Thus as the black contours in the range $(-2,1)$ approach the real axis,
$$
\log\zeta(s_2)=2\pi i+\log\zeta(s_1)
$$ 
and therefore
$$
\text{Im}\log\zeta(s_2)\to 2\pi+\text{Im}\log\zeta(s_1).
$$

Now, the real part of $\zeta(s)$ in the range $(-2,1)$ is negative and the imaginary part of $\zeta(s)$ is negative in the range $(-2+i\epsilon,2+i\epsilon)$ along the upper black contour (gray contour in Figure \ref{fig:realaxiscloseup}).  Therefore $\arg\zeta(s)$ along the upper black contour in the range $-2<\sigma<1$ must approach $-\pi$ since this branch of $\arg$ is continuous and therefore by the expression above, $\arg\zeta(s)$ along the lower contour is $\pi$ over the same range.

Now, in the vicinity of each trivial zero, 
\begin{equation}
\frac{\zeta'(s)}{\zeta(s)}=\frac{1}{s+2n}+\frac{\phi'(s)}{\phi(s)}
\label{eqn0035}
\end{equation}
with the residue of this expression clearly being one.  Referring to Figure (\ref{fig:closeupat2}), consider the points $s_3$ and $s_4$.  Then
\begin{equation*}
\log\zeta(s_4)=\lim_{r\to 0}\left\{\int_{\sigma_1+ir}^{-2+r+ir}\frac{\zeta'(s)}{\zeta(s)}ds+\mathop\int\limits_{\Gamma_r}\frac{\zeta'(s)}{\zeta(s)}ds+\int_{-2-r+ir}^{\sigma_2+ir}\frac{\zeta'(s)}{\zeta(s)}ds\right\}+\log\zeta(s_3)
\end{equation*}
and again because of the symmetric form of the limit,
$$\lim_{r\to 0}\mathop\int\limits_{\Gamma_r}\frac{\zeta'(s)}{\zeta(s)}ds=\frac{\pi}{2} i,$$
leaving for the integrals:
\begin{equation}
\lim_{r\to 0}\left\{\int_{\sigma_1+ir}^{-2+r+ir}\frac{\zeta'(s)}{\zeta(s)}ds+\int_{-2-r+ir}^{\sigma_2+ir}\frac{\zeta'(s)}{\zeta(s)}ds\right\}.
\label{eqn0037}
\end{equation}
Substituting  \eqref{eqn0035} into  \eqref{eqn0037}, we obtain after noting the analytic parts cancel:
\begin{equation*}
\begin{aligned}
\lim_{r\to 0}\left\{\int_{\sigma_1+ir}^{-2+r+ir}\frac{1}{s+2}ds+\int_{-2-r+ir}^{\sigma_2+ir}\frac{1}{s+2}ds\right\}&\\
&\hspace{-150pt}=\lim_{r\to 0}\Big\{\Log(r+ir)-\Log(\sigma_1+2+ir)-\Log(-r+ir)+\Log(2+\sigma_2+ir)\Big\}\\
&\hspace{-150pt}=\frac{\pi}{4}i+0-\frac{3\pi}{4}i+\pi i -|\sigma_1+2|+|\sigma_2+2|\\
&\hspace{-150pt}=\frac{\pi}{2}i+f(\sigma_1,\sigma_2)
\end{aligned}
\end{equation*}
since $\sigma_2<-2<\sigma_1$.  Then
$$
\log\zeta(s_4)=\pi i+f(\sigma_1,\sigma_2)+\log\zeta(s_3),
$$
and
$$\text{Im}\log\zeta(s_4)=\pi+\text{Im}\log\zeta(s_3).
$$
Therefore the imaginary component of $\log\zeta(s)$ along the upper black contour as it approaches the real axis in the range $(-4,-2)$ has the limit zero since we have shown above that $\text{Im}\log\zeta(s)$ for $-2<\sigma<1$ along the upper black contour has a limit of $-\pi$.  This same approach can then be applied to the remaining components of the black contour to arrive at the geometry in Figure (\ref{fig:realaxiscloseup}).

\begin{figure}
\centering
%%----start of first subfigure----
\subfloat[Common zero]{
\label{fig110:subfig:a} %% label for first subfigure
\includegraphics[scale=.7]{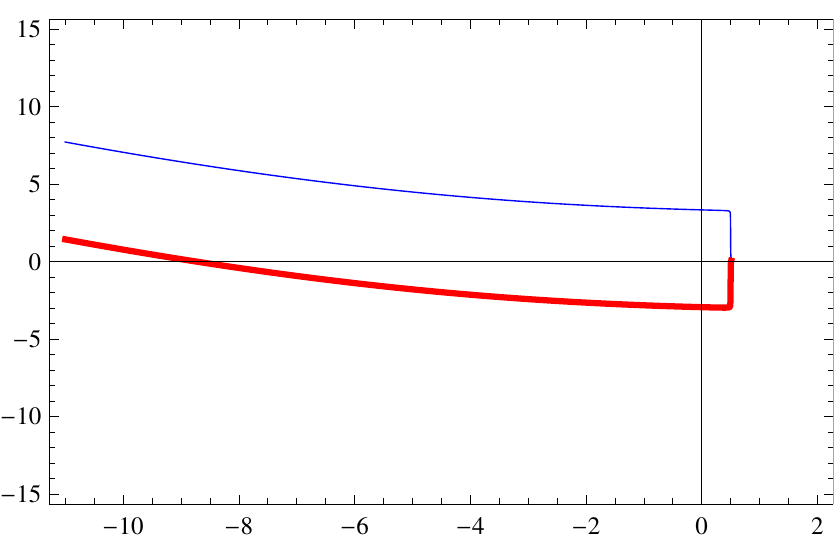}}
\hspace{1in}
%%----start of second subfigure----
\subfloat[Rogue zeros]{
\label{fig110:subfig:b} %% label for second subfigure
\includegraphics[scale=.7]{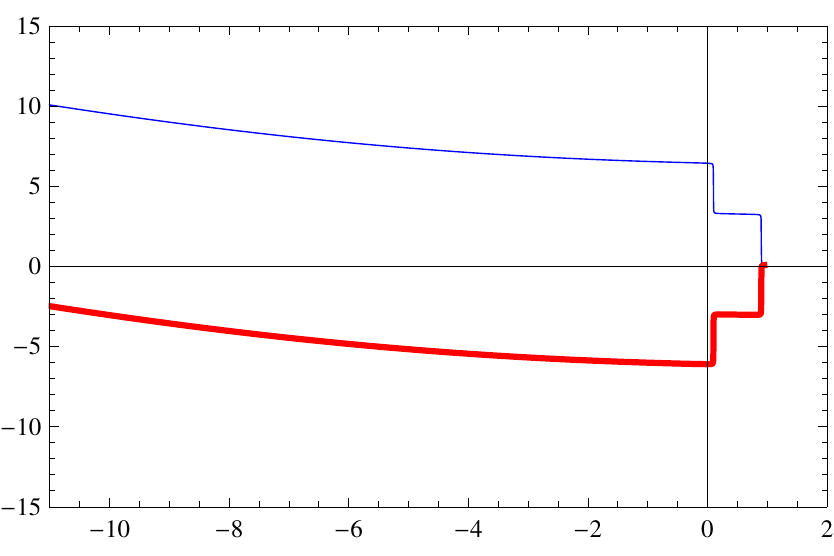}}
\caption{$\arg\zeta(s)$ along red and blue contours of $B(m)$}
\label{fig110:subfig} %% label for entire figure
\end{figure}

\subsection{Geometry of $\log\zeta(s)$ along branch cuts at each non-trivial zero}
\label{subsection32}
$\arg\zeta(s)$ in a region around the branch cut at $\rho_1$ is shown in Figure  (\ref{fig110:subfig:a}). 
Empirical results indicate the Red contour (over the imaginary sheet of $\log\zeta(s)$) is $-2\pi$ that of the Blue contour.  A similar argument to that in Section \ref{section31} can be used to explain this geometry: For a point $s_1=\sigma_1+i(\gamma_n+r)$ on the blue contour and $s_2=\sigma_1+i(\gamma_n-r)$ on the red contour, 
\begin{equation*}
\log\zeta(s_2)=\lim_{r\to 0}\left\{\int_{\sigma_1+i(\gamma_n+r)}^{1/2-r+i(\gamma_n+r)}\frac{\zeta'(s)}{\zeta(s)}ds+\mathop\int\limits_{\Gamma_r}\frac{\zeta'(s)}{\zeta(s)}ds+\int_{1/2-r+i(\gamma_n-r)}^{\sigma_1+i(\gamma_n-r)}\frac{\zeta'(s)}{\zeta(s)}ds\right\}+\log\zeta(s_1).
\end{equation*}
And because of symmetry,
$$
\lim_{r\to 0}\mathop\int\limits_{\Gamma_r}\frac{\zeta'(s)}{\zeta(s)}ds=-\frac{3\pi}{2}i.
$$
For a deleted neighborhood surrounding each trivial zero on the critical line, called a ``common zero'' in this paper, we have
\begin{equation*}
\frac{\zeta'(s)}{\zeta(s)}=\frac{1}{s-\rho_n}+\frac{\phi'(s)}{\phi(s)}.
\end{equation*}
Making the same type of substitution as above gives
\begin{equation*}
\begin{aligned}
\lim_{r\to 0}\left\{\int_{\sigma_1+i(\gamma_n+r)}^{1/2-r+i(\gamma_n+r)}\frac{1}{s-\rho_n}ds+\int_{1/2-r+i(\gamma_n-r)}^{\sigma_1+i(\gamma_n-r)}\frac{1}{s-\rho_n}ds\right\}&\\
&\hspace{-200pt}=\lim_{r\to 0}\Big\{\Log(-r+ir)-\Log(\sigma_1-1/2+ir)-\Log(-r-ir)+\Log(\sigma_1-1/2-ir)\Big\} \\
&\hspace{-200pt}=\frac{3\pi}{4}i-\pi i+\frac{3\pi}{4}i-\pi i=-\frac{\pi}{2}i.
\end{aligned}
\end{equation*}
This leaves
$$
\log\zeta(s_2)\to -2\pi i+\log\zeta(s_1)
$$
and thus
$$
\text{Im}\log\zeta(s_2)\to -2\pi+\text{Im}\log\zeta(s_1)
$$
which is the value observed empirically.

We must also consider the possibility of zeros off the critical line.  These are called "rogue zeros" in this paper.  How would these effect the branch cuts?  We can artificially introduce such zeros by considering the function $\zeta(s)(s-r_1)(s-r_2)$ and numerically solve (\ref{inteqn1}) using this function.  The plot in Figure (\ref{fig110:subfig:b}) shows $\arg\left[\zeta(s)(s-r_1)(s-r_2)\right]$ at this branch cut with $r_1=.1+i\gamma$ and $r_2=.9+i\gamma$.  Numerical results suggest the difference in argument from $-\sigma$ to the first zero is $4\pi$ and between the zeros is $2\pi$.  This same type of geometry would occur with actual rogue zeros and using the same argument above, this $4\pi$ and $2\pi$ difference is easily explained.

\section{Integral setup}
We now consider
$$
\mathop\oint\limits_{B(m)}\frac{\log\zeta(s)}{s}x^s ds
$$
where $B(m)$ includes the possibility of zeros off the critical line with the horizontal gray legs of the contour passing over certain lines with imaginary components equal to $\pm T_m$. Since $B(m)$ is closed over a function analytic throughout the region of integration,
$$
\mathop\oint\limits_{B(m)}\frac{\log\zeta(s)}{s}x^s ds=0,
$$
and therefore
\begin{equation*}
\begin{aligned}
\underset{B(m)}{\oint }\log\zeta (s)\frac{x^s}{s} ds&=\underset{\text{Purple}}{\int }\log  \zeta (s)\frac{x^s}{s} ds+\underset{\text{Green}}{\int }\log  \zeta (s)\frac{x^s}{s} ds\\
&+\underset{\text{Cyan} }{\int }\log  \zeta (s)\frac{x^s}{s} ds,\underset{\text{Orange} }{\int }\log  \zeta (s)\frac{x^s}{s} ds+\underset{\text{Black} }{\int }\log  \zeta (s)\frac{x^s}{s} ds \\
&+\underset{\text{Brown}}{\int }\log  \zeta (s)\frac{x^s}{s} ds+\underset{\text{red}+\text{blue}}{\int }\log  \zeta (s)\frac{x^s}{s} ds+\underset{\text{Gray} }{\int }\log  \zeta (s)\frac{x^s}{s} ds \\&=0.
\end{aligned}
\end{equation*}

Now $\log\zeta(s)=\text{Log}\,\zeta(s)$ over the purple contour and therefore we can replace $\log\zeta(s)$ with $\text{Log}\,\zeta(s)$ in that integral. Multiplying by $\frac{1}{2\pi i}$,  and taking the limit as $m\to\infty$,
\begin{equation}
\begin{aligned}
\frac{1}{2\pi  i}\underset{a-\text{i$\infty $}}{\overset{a+\text{i$\infty $}}{\int }}\text{Log}\,\zeta (s)\frac{x^s}{s} ds
&=-\frac{1}{2\pi  i}\left\{\underset{\text{Green}}{\int }\log  \zeta (s)\frac{x^s}{s} ds+\underset{\text{Cyan}}{\int }\log  \zeta (s)\frac{x^s}{s} ds\right.\\
&\hspace{25pt}+\underset{\text{Orange}}{\int }\log  \zeta (s)\frac{x^s}{s} ds+\underset{\text{Black}}{\int }\log  \zeta (s)\frac{x^s}{s} ds+\underset{\text{Brown}}{\int }\log  \zeta (s)\frac{x^s}{s} ds \\
&\hspace{25pt}\left.+\underset{\text{Red}+\text{Blue}}{\int }\log  \zeta (s)\frac{x^s}{s} ds+\underset{\text{Gray}}{\int }\log  \zeta (s)\frac{x^s}{s} ds\right\}.
\end{aligned}
\label{eqn100}
\end{equation}
Where in  \eqref{eqn100}, the appropriate limits are taken for each integral.  

We now go on to evaluate the various color-coded components of this expression.  

\section{Evaluation of the integral over B(m)}
%
% contour over green circles
%
\subsection{Integral over the green contours}

The green contours encircle each non-trivial zero of the zeta function.  We allow the diameter of these contours to approach zero and take the following limit:

$$\lim_{\epsilon \rightarrow 0} \underset{\text{Green}}{\int }\frac{\log  \zeta (s)}{s}x^s ds;\quad s=\rho _n+\epsilon  e^{\text{it}}$$
where $\left\{\rho _n\right\}$ represents the set of non-trivial zeros of the zeta function.  Writing this in terms of s=$\rho_n+\epsilon e^{\text{i$\theta $}}$, we have for the common zeros:
$$
\underset{\text{Green(c)}}{\int }\frac{\log  \zeta (s)}{s}x^s ds= \int _{-\pi }^{\pi }\frac{\log  \zeta \left(\rho _n+\epsilon  e^{i \theta }\right)}{\rho _n+\epsilon  e^{i \theta }}\epsilon  i e^{i \theta } x^s d\theta
$$
and
$$
\underset{\text{Green(r)}}{\int }\frac{\log  \zeta (s)}{s}x^s ds= \int _{\pi }^{0 }\frac{\log  \zeta \left(\rho _n+\epsilon  e^{i \theta }\right)}{\rho _n+\epsilon  e^{i \theta }}\epsilon  i e^{i \theta } x^s d\theta+ \int _{0 }^{-\pi }\frac{\log  \zeta \left(\rho _n+\epsilon  e^{i \theta }\right)}{\rho _n+\epsilon  e^{i \theta }}\epsilon  i e^{i \theta } x^s d\theta
$$
for rogue zeros with $\sigma_1<\sigma_2$.

The $x^s$ term is bounded by $\left|x^s\right|$ and note as $\epsilon\to 0$, the denominator tends to $\rho _n$ with the numerator tending to
$$
\epsilon k\log  \zeta \left(\rho +\epsilon  e^{\text{it}}\right)
$$
with $k$ being an upper bound on the remaining terms. The zeta function has simple zeros and therefore
$$
\lim_{\epsilon\to 0} \epsilon k \log\zeta(\rho+\epsilon e^{it})\sim \lim_{\epsilon\to 0} \epsilon k \log (\epsilon) 
$$
which tends to zero.  This gives
\begin{equation*}
\lim_{\epsilon\to0}\mathop\int\limits_{\text{Green}} \frac{log\zeta(s)}{s} x^s ds=0
\end{equation*}
independent of the placement of non-trivial zeros.

%
% yellow contour
%
%
% yellow contour
%
\subsection{Integral over the cyan contours}

The cyan contours encircle the trivial zeros and yield the following integrals:

$$\lim_{\epsilon \rightarrow 0} \underset{\text{Cyan}}{\int }\frac{\log  \zeta (s)}{s}x^s ds,\quad s=-2n+\epsilon  e^{i\theta}.$$
Similar to the integrals over the green indentations, we express these integrals in terms of s=$-2n+\epsilon e^{i\theta}$ and since these zeros are simple, we obtain a similar limiting process:
$$\lim_{\epsilon \rightarrow 0}  \epsilon k \log  \zeta \left(-2n+\epsilon  e^{i\theta}\right)\sim\lim_{\epsilon\to0} \epsilon k \log \epsilon,\quad n=-2,-4,\cdots $$
and therefore

\begin{equation*}
\lim_{\epsilon\to 0}\underset{\text{Cyan}}{\int }\frac{\log  \zeta (s)}{s}x^s\text{ds}=0.
\end{equation*}
%
% orange contour
%
% orange contour
\subsection{Integral over the orange contour}
The integral around the orange contour at s=1 yields
$$
\lim_{\epsilon \rightarrow 0} \underset{\text{Orange}}{\int }\frac{\log  \zeta (s)}{s}x^s ds,\quad s=1+\epsilon  e^{i\theta}.$$
Similar to the integrals over the green indentations, we express the integral in terms of s=$1+\epsilon e^{i\theta}$ and arrive at a similar limiting process:
$$
\lim_{\epsilon \rightarrow 0}  \epsilon k \log  \zeta \left(1+\epsilon  e^{i\theta}\right).
$$
The pole at s=1 is a simple pole so that this limit approaches the form
$$
\lim_{\epsilon \rightarrow 0}  \epsilon  \log  \zeta \left(1+\epsilon  e^{\text{it}}\right)\sim \lim_{\epsilon \rightarrow 0}  \epsilon  \log \text{  }\frac{1}{\epsilon }=0.
$$
Therefore

\begin{equation*}
\lim_{\epsilon\to0}\underset{\text{Orange}}{\int }\frac{\log  \zeta (s)}{s}x^s ds=0.
\end{equation*}
%

%
% **** brown contour around the origin
%
\subsection{Integral over the brown contour around the origin}

We wish to consider the following limit:
$$
\lim_{\epsilon \rightarrow 0} \underset{\text{Brown}}{\int }\frac{\log  \zeta (s)}{s}x^s ds;\quad s=\epsilon  e^{\text{it}}.
$$
Referring to Figure (\ref{fig:file1}), we can write this in terms of two integrals: 
$$
\underset{\text{Brown}}{\int }\frac{\log  \zeta (s)}{s}x^s ds= \int _{\pi }^0\frac{\log  \zeta \left(\epsilon  e^{\text{it}}\right)}{\epsilon  e^{\text{it}}}\epsilon  i e^{\text{it}} x^{\gamma }dt+\int _{0}^{-\pi}\frac{\log  \zeta \left(\epsilon  e^{\text{it}}\right)}{\epsilon  e^{\text{it}}}\epsilon  i e^{\text{it}} x^{\gamma }dt.
$$

Note each integral traverses surfaces separated by a branch cut with $\Delta\arg\zeta(s)=2\pi$ as described in Section \ref{section31}.  As discussed in that section, $\arg\zeta$ over the upper contour is $-\pi$ and over the lower contour is $\pi$ as shown in Figure (\ref{fig:realaxiscloseup}).   We then have for the upper contour:
\begin{equation*}
\begin{aligned}
\lim_{\epsilon \rightarrow 0} \int _{\pi }^0\frac{\log  \zeta \left(\epsilon  e^{\text{it}}\right)}{\epsilon  e^{\text{it}}}\epsilon  i e^{\text{it}} x^{\gamma }dt
&=\lim_{\epsilon \rightarrow 0} i\int _{\pi }^0\log  \zeta \left(\epsilon  e^{\text{it}}\right) x^{\gamma }dt\\
&=\lim_{\epsilon \rightarrow 0} i\int _{\pi }^0\Big(\log \left|\zeta \left(\epsilon  e^{\text{it}}\right)\right|+i \arg  \zeta \left(\epsilon  e^{\text{it}}\right)\Big)x^{\gamma} \text{dt}
\end{aligned}
\end{equation*}
with $\gamma=\epsilon e^{it}$.
As $\epsilon \to 0$, the first term of the last integral tends to $\log|\zeta(0)|$ and the imaginary component of $\log\zeta $(s) over the contour near zero tends to $\arg\zeta(0)=-\pi$.  That is,
$$
\lim_{\epsilon \rightarrow 0}  \arg  \zeta \left(\epsilon  e^{\text{it}}\right)\rightarrow -\pi,\quad 0< t< \pi. 
$$

Using these values we have
\begin{equation*}
\begin{aligned}
\lim_{\epsilon \rightarrow 0} i\int _{\pi }^0\log  \zeta \left(\epsilon  e^{\text{it}}\right) x^{\epsilon  e^{\text{it}}}dt&=-\pi  i(\log |\zeta (0)|-\pi  i)\\ &=-\pi  i(-\log (2)-\pi  i)\\
&=-\pi ^2+\pi  i \log  2.
\end{aligned}
\end{equation*}

In the same way, $\arg\zeta(s)$ over the lower contour approaches $\pi $ according to Section \ref{section31}, and we obtain for that result:

\begin{equation*}
\begin{aligned}
\lim_{\epsilon \rightarrow 0} \int _0^{-\pi }\log  \zeta \left(\epsilon  e^{\text{it}}\right) x^{\epsilon  e^{\text{it}}}dt&=-\pi  i(\log |\zeta (0)|+\pi  i)\\&=-\pi  i(-\log (2)+\pi  i)\\&=\pi ^2+\pi  i \log  2.
\end{aligned}
\end{equation*}

Therefore the sum of the integrals over the central brown contours as the diameters tend to zero is

\begin{equation*}
\lim_{\epsilon\to0} \underset{\text{Brown}}{\int }\frac{\log  \zeta (s)}{s}x^s ds=2\pi  i \log (2).
\end{equation*}

\subsection{Integrals over the blue and red contours}
We consider two cases:  Case I takes the zeros on the critical line and Case II considers zeros off the critical line.  

\subsubsection{Zeros on the critical line}
In the case of zeros on the critical line,  the value of log $\zeta $(s) on the blue contour is $\log|\zeta(s)|+i \arg \zeta(s)$ and the value on the red contour is $\log|\zeta(s)|+i( \arg\zeta(s)-2\pi)$ as discussed in Section \ref{subsection32}.  In the limit as $s\to \sigma+i\gamma_n$, the integrals over these contours as $B(m)$ extends to infinity are

\begin{equation*}
\begin{aligned}
\underset{\text{Blue(n)}}{\int }\frac{\log  \zeta (s)}{s}x^s ds+\underset{\text{Red(n)}}{\int }\frac{\log  \zeta (s)}{s}x^s ds&\\
&\hspace{-125pt}=\lim_{\epsilon\to 0}\left\{\int _{-\infty }^{1/2-\epsilon}\frac{\log |\zeta (s)|+i \arg  \zeta (s)}{s}x^s d\sigma-\int _{-\infty }^{1/2-\epsilon}\frac{\log |\zeta (s)|+i(\arg  \zeta (s)-2\pi )}{s}x^s d\sigma\right\} \\
&\hspace{-125pt}=2\text{$\pi $i}\int_{-\infty }^{1/2} \frac{e^{\left(\sigma+i \gamma _n\right) \log (x)}}{\sigma+i\gamma_n} \, d\sigma.
\end{aligned}
\end{equation*}

Letting $u=\sigma+\text{i$\gamma $}_n$, then

$$2\text{$\pi $i}\int_{-\infty +\text{i$\gamma_n$}}^{1/2+\text{i$\gamma_n $}} \frac{e^{u \log (x)}}{u} du,
$$
and with $v=u\log(x)$ gives
$$\underset{\text{Blue(n)}}{\int }\frac{\log  \zeta (s)}{s}x^s ds+\underset{\text{Red(n)}}{\int }\frac{\log  \zeta (s)}{s}x^s ds=2\pi  i\underset{-\infty +i \gamma_n  \log (x)}{\overset{1/2 \log (x)+i \gamma_n  \log (x)}{\int }}\frac{e^v}{v} dv=2\pi i\left[Li(x^{\rho})-\pi i\right].
$$
For the conjugate contours at $\rho_{-n}$,
$$\underset{\text{Blue(-n)}}{\int }\frac{\log  \zeta (s)}{s}x^s ds+\underset{\text{Red(-n)}}{\int }\frac{\log  \zeta (s)}{s}x^s ds=2\pi  i\underset{-\infty -i \gamma_n  \log (x)}{\overset{1/2 \log (x)-i \gamma_n  \log (x)}{\int }}\frac{e^v}{v} dv=2\pi i\left[Li(x^{1-\rho})+\pi i\right]
$$
with the last integrals expressed in terms of the logarithmic integral as per   \eqref{ingham001}.

\subsubsection{Zeros off the critical line}

In this case we have sets of zeros of the form

$$
\left\{\rho _n\right\}=\left\{\sigma _1\pm \text{i$\gamma $}_n,\sigma _2\pm \text{i$\gamma $}_n\right\};\text{     }0<\sigma _1<1/2<\sigma _2<1.
$$
The presence of these zeros will cause the values between the red and blue contours to differ by $4\pi$ between the range $-\infty$ to $\sigma_1$ and by $2\pi i$ between the range $(\sigma_1,\sigma_2)$ as was discussed in Section \ref{subsection32}.  We obtain for each set of zeros $\sigma_1+i\gamma_n, \sigma_2+i\gamma_n$ and with $s=\sigma+i\gamma_n$,

\begin{equation*}
\begin{aligned}
\underset{\text{Blue(n)}}{\int }\frac{\log  \zeta (s)}{s}x^s ds+\underset{\text{Red(n)}}{\int }\frac{\log  \zeta (s)}{s}x^s ds\\
&\hspace{-150pt}=\lim_{\epsilon\to 0}\left\{\int _{-\infty }^{\sigma _1-\epsilon}\frac{\log |\zeta (s)|+i \arg  \zeta (s)}{s}x^s d\sigma
-\int _{-\infty }^{\sigma _1-\epsilon}\frac{\log |\zeta (s)|+i(\arg  \zeta (s)-4\pi )}{s}x^s d\sigma\right. \\
&\hspace{-130pt}+\left.\int _{\sigma _1+\epsilon}^{\sigma _2-\epsilon}\frac{\log |\zeta (s)|+i (\arg  \zeta (s)-\pi )}{s}x^s d\sigma-\int _{\sigma _1+\epsilon}^{\sigma _2-\epsilon}\frac{\log |\zeta (s)|+i(\arg  \zeta (s)-3\pi )}{s}x^s d\sigma\right\}.
\end{aligned}
\end{equation*}
This then becomes

$$
\begin{aligned}
4\pi  i\int_{-\infty }^{\sigma _1} \frac{x^s}{s} \, d\sigma+2\pi  i \int_{\sigma _1}^{\sigma _2} \frac{x^s}{s} \, d\sigma
&=2\text{$\pi $i} \int_{-\infty }^{\sigma _1} \frac{x^s}{s} \, d\sigma+2\text{$\pi $i} \int_{-\infty }^{\sigma _1} \frac{x^s}{s} \, d\sigma+2\pi  i \int_{\sigma _1}^{\sigma _2} \frac{x^s}{s} \, d\sigma \\
&=2\text{$\pi $i} \int_{-\infty }^{\sigma _1} \frac{x^s}{s} \, d\sigma+2\text{$\pi $i} \int _{-\infty }^{\sigma _2}\frac{x^s}{s}d\sigma;\quad s=\sigma+i\gamma_n
\end{aligned}
$$
with a similar case for the pair with negative imaginary parts.  These expressions are again  equal to $Li(x^\rho)-\pi i$. That is,
$$
\int_{-\infty}^{\sigma_n}\frac{x^s}{s} d\sigma=\underset{-\infty +i \gamma_n  \log (x)}{\overset{\sigma_n \log (x)+i \gamma_n  \log (x)}{\int }}\frac{e^v}{v}\text{dv}= Li(x^{\rho_n})-\pi i;\quad s=\sigma+i\gamma_n, \quad \rho_n=\sigma_n+i\gamma_n.$$
The set of conjugate zeros gives a similar expression in terms of $Li(x^{1-\rho_n})+\pi i$. Therefore we can write for the set of common and rogue zeros,
\begin{equation*}
\underset{\substack{\text{Red(n)}\\ \text{Blue(n)}}}{\int }\frac{\log  \zeta (s)}{s}x^s ds=\begin{cases}\displaystyle2\pi i\int_{-\infty }^{1/2} \frac{x^s}{s} \, d\sigma &\mbox{\text{(common zeros)}}\vspace{10pt} \\ 
\displaystyle 2 \pi i \int_{-\infty}^{\sigma_1}\frac{x^s}{s}d\sigma+2\pi i \int_{-\infty}^{\sigma_2}\frac{x^s}{s}d\sigma & \mbox{\text{(rogue zeros)}}.
\end{cases}
\end{equation*}
So that for the set of all non-trivial zeros $\{\rho\}$: 

\begin{equation*}
\sum_{\rho}\hspace{10pt}\mathop\int\limits_{\substack{\text{Red(n)} \\ \text{Blue(n)}}}\frac{\log\zeta(s)}{s}x^s 
ds=\sum_{\text{Im}(\rho_n)>0}\Big[Li(x^{\rho_n})+Li(x^{1-\rho_n})\Big].
\end{equation*}

\subsection{Integral over the black contours between $-\infty$ and -4}

The integrals over the black contours are determined by the change in argument over these contours described in Section \ref{section31}. Those over the interval (-6,-4) are 
\begin{equation*}
\begin{aligned}
\mathop\int\limits_{\substack{\text{upper} \\ \text{contour}}} \frac{\log\zeta(s)}{s}x^s ds+\mathop\int\limits_{\substack{\text{lower} \\ \text{contour}}} \frac{\log\zeta(s)}{s}x^s ds & \\
&\hspace{-150pt}=
\lim_{\epsilon\to 0}\left\{\int_{-6+\epsilon}^{-4-\epsilon}\frac{\log |\zeta (s)|+i \arg  \zeta (s)}{s}x^s ds\right. +\left.\int_{-4-\epsilon}^{-6+\epsilon}\frac{\log |\zeta (s)|+i (\arg  \zeta (s)-2\pi )}{s}x^s ds\right\} \\
&\hspace{-150pt}=2\pi  i \int_{-6}^{-4} \frac{x^s}{s} \, ds=2\pi  i [\Gamma (0,6 \log  x)-\Gamma (0,4 \log  x)],
\end{aligned}
\end{equation*}
with $\Gamma(0,x)$ being the incomplete gamma function. Over the interval $(-6,-8)$, these integrals are 
\begin{equation*}
\begin{aligned}
\lim_{\epsilon\to 0}\left\{\int_{-8+\epsilon}^{-6-\epsilon}\frac{\log |\zeta (s)|+i \arg  \zeta (s)}{s}x^s ds\right. &+\left.\int_{-6-\epsilon}^{-8+\epsilon}\frac{\log |\zeta (s)|+i (\arg  \zeta (s)-4\pi )}{s}x^s ds\right\} \\
&\hspace{-50pt}=4\pi  i \int_{-8}^{-6} \frac{x^s}{s} \, ds=4\pi  i [\Gamma (0,8 \log  x)-\Gamma (0,6 \log  x)].
\end{aligned}
\end{equation*}
And for the interval $(-8,-10)$,
\begin{equation*}
\begin{aligned}
\lim_{\epsilon\to 0}\left\{\int _{-10+\epsilon}^{-8-\epsilon}\frac{\log |\zeta (s)|+i \arg  \zeta (s)}{s}x^sds\right.&\left.+\int _{-8-\epsilon}^{-10+\epsilon}\frac{\log |\zeta (s)|+i (\arg  \zeta (s)-6\pi )}{s}x^sds\right\}\\
&\hspace{-50pt}=6\pi  i \int_{-10}^{-8} \frac{x^s}{s} \, ds=6\pi  i [\Gamma (0,10 \log  x)-\Gamma (0,8 \log  x)].
\end{aligned}
\end{equation*}
And in general,
\begin{equation*}
\begin{aligned}
\lim_{\epsilon\to 0}\left\{\int _{-2n+\epsilon}^{-2n+2-\epsilon}\frac{\log |\zeta (s)|+i \arg  \zeta (s)}{s}x^sds\right.&\left.+\int _{-2n+2-\epsilon}^{-2n+\epsilon}\frac{\log |\zeta (s)|+i (\arg  \zeta (s)-(2n-4)\pi )}{s}x^sds\right\}\\
&\hspace{-50pt}=2n \pi  i \int_{-10}^{-8} \frac{x^s}{s} \, ds=2\text{n$\pi $} i [\Gamma (0,2n \log  x)-\Gamma (0,(2n-2) \log  x)].
\end{aligned}
\end{equation*}

Adding the first few terms of this set gives
\begin{equation}
\begin{aligned}
2\pi  i \Gamma (0,6 \log  x)-2\pi  i \Gamma (0,4 \log  x)&+4\pi  i \Gamma (0,8 \log  x)-4\pi  i \Gamma (0,6 \ln  x)\\
&+6\pi  i \Gamma (0,10 \log  x)-6\pi  i \Gamma (0,8 \log  x)+\cdots\\ 
&+2\text{n$\pi $} i \Gamma (0,2n \log  x)-2n\pi  i \Gamma (0,2n-2 \log  x)+\cdots. 
\end{aligned}
\end{equation}
Summing all the terms, we have for the black contours in the interval $(-\infty,-4)$:
\begin{equation}
\mathop\int\limits_{\text{Black}(-\infty,-4)}\frac{\log  \zeta (s)}{s}x^s ds=-2\pi  i \sum _{n=2}^{\infty }  \Gamma [0,2n \log (x)],
\end{equation}
with the sum converging due to the exponential decay of the gamma integrals.

\subsection{Integral over the black contour from -2 to 1}
As discussed in Section \ref{section31}, in the range (-2,1), the imaginary surface over the lower contour is $\pi$ and over the upper contour, $-\pi$.  Therefore, log $\zeta $(s) differs by 2$\pi $i over these two contours.  However, the integrand is singular at the origin and thus we take the principal value for these integrals.  We then have
\begin{equation*}
\begin{aligned}
\mathop\int\limits_{\text{Black}(-2,1)}\frac{\log  \zeta (s)}{s}x^s ds&=\mathop\int\limits_{\substack{\text{upper}\\ \text{contour}}}\frac{\log  \zeta (s)}{s}x^s ds+\mathop\int\limits_{\substack{\text{lower}\\ \text{contour}}}\frac{\log  \zeta (s)}{s}x^s ds \\
&\hspace{-75pt}=P.V.\left\{\lim_{\epsilon\to 0}\left(\int_{-2+\epsilon}^{1-\epsilon}\frac{\log |\zeta (s)|+i \arg  \zeta (s)}{s}x^s ds-\int_{-2+\epsilon}^{1-\epsilon}\frac{\log |\zeta (s)|+i (\arg  \zeta (s)+2\pi )}{s}x^s ds\right)\right\} \\
&\hspace{-75pt}=-2\pi i\; P.V. \int _{-2}^1\frac{x^s}{s} ds.
\end{aligned}
\end{equation*}
\subsection{Integral over the gray contours}
We wish to prove the following:
\begin{equation*}
\lim_{m\to\infty}\mathop\int\limits_{\text{Gray}}\frac{\log\zeta(s)}{s}x^s ds=0
\end{equation*}
where the gray boundary is incremented through certain lines crossing the critical strip of the zeta function .  Proofs of the following theorems can be found in Ingham \cite{Ingham}.

\begin{theorem}
There exist a sequence of numbers, $T_2,T_3,\cdots$, such that:
\begin{equation*}
\begin{aligned}
m&<T_m<m+1; \quad (m=2,3,\cdots) \\
\left|\frac{\zeta'(s)}{\zeta(s)}\right|&<A\log^2 T_m,\quad s=\sigma+iT_m,\quad -1\leq\sigma\leq 2
\end{aligned}
\end{equation*}
\label{thm001}
for some positive constant $A$.
\end{theorem}

These are the lines at $\pm T_m$ of Figure (\ref{fig:file1}).

\begin{theorem}
In the region obtained by removing from the half-plane $\sigma\leq -1$, the interiors of a set of circles surrounding each trivial zero with radius $1/2$, we have:
\begin{equation*}
\left|\frac{\zeta'(s)}{\zeta(s)}\right|<A \log(|s|+1)<A\log^2 T_m;\quad \sigma<-1
\end{equation*}
\label{thm002}
\end{theorem}

And for $\sigma>2$:
\begin{equation*}
\left|\frac{\zeta'(s)}{\zeta(s)}\right|=\left|\sum_{n=1}^{\infty}\frac{\Lambda(n)}{n^s}\right|<\sum_{n=1}^{\infty}\frac{\text{Log}(n)}{n^{2}}.
\end{equation*}
\subsubsection{Integrals over horizontal legs of gray contour}

Theorem $(\ref{thm001})$ implies that between each integer $m$ and $m+1$, there exists a horizontal contour passing between the zeta zeros over which the order of $\frac{\zeta'(s)}{\zeta(s)}$ does not exceed $A\log^2 T_m$.  Theorem $(\ref{thm002})$ states this order does not exceed $A\log(|s|+1)$ when $\text{Re}(s)<-1$.
\begin{figure}
	\centering
		\includegraphics{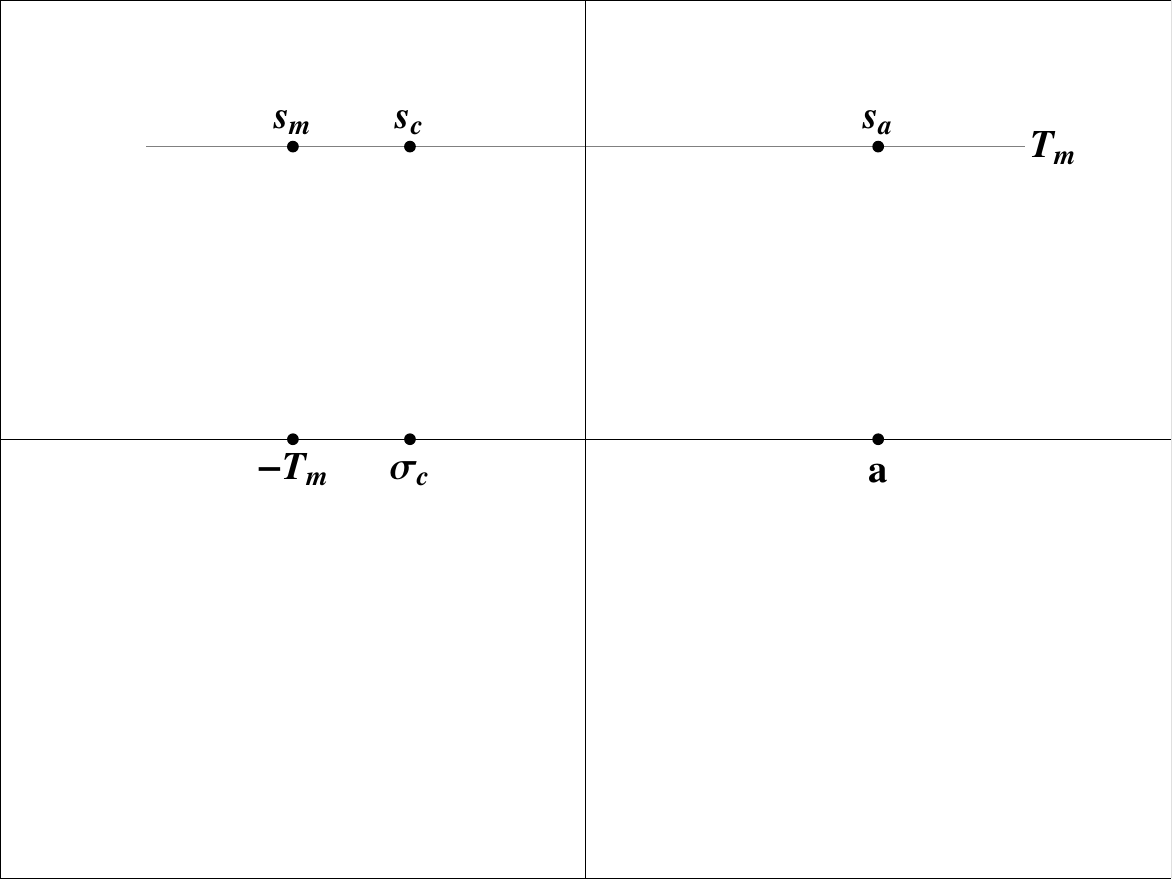}
	\caption{$T_m$ contour (for $T_m$ sufficiently large)}
	\label{fig:tmcontour}
\end{figure}

Figure(\ref{fig:tmcontour}) shows the top gray contour, $\text{Gray($T_m$)}$, over a finite path of length $|-T_m-a|$.  Then

\begin{equation}
\mathop\int\limits_{\text{Gray($T_m$)}}\frac{\log\zeta(s)}{s}x^s ds=\int_{s_a}^{s_c}\frac{\log\zeta(s)}{s}x^s ds+\int_{s_c}^{s_m}\frac{\log\zeta(s)}{s}x^s ds;\quad \text{Re}(s_m)=-T_m
\label{eqn512}
\end{equation}
where $\sigma_c=-\frac{3}{\log(x)}\log T_m=C_2 \log T_m$.  That is, for $\sigma$ sufficiently negative and $T_m$ sufficiently large, we have the inequality

\begin{equation*}
|x^s|<(T_m)^{-3},\quad \sigma<\sigma_c\;\text{and}\,-T_m<\sigma_c.
\end{equation*}

The value of $\log\zeta(s)$ over this contour given by  \eqref{inteqn1} is
\begin{equation*}
\log\zeta(s)=\int_{s_a}^s\frac{\zeta'(s)}{\zeta(s)}ds+\log\zeta(s_a).
\end{equation*} 
And therefore
\begin{equation}
\begin{aligned}
\Bigg.|\log\zeta(s)|\Bigg|_{(s_c,s_a)}&\leq\int_{s_a}^{s_c}\left|\frac{\zeta'(s)}{\zeta(s)}\right|ds+K\\
&\leq |s_a-s_c|(A\log^2 T_m+K) \\
&\leq |\sigma_c-a|(A\log^2 T_m+K) \\
&\leq |C_2 \log T_m-a|(A\log^2 T_m+K),
\end{aligned}
\label{eqnlog1}
\end{equation}
where $A$ is a positive constant and $K=\log\sum\frac{\log(n)}{n^{1+\delta}}$ for $1+\delta<a$.  Substituting  \eqref{eqnlog1} into  \eqref{eqn512}, we can write for the first integral:
\begin{equation}
\begin{aligned}
\left|\int_{s_a}^{s_c}\frac{\log\zeta(s)}{s}x^s ds\right|&\leq\int_{s_a}^{s_c}\frac{|C_2 \log T_m-a|(A\log^2 T_m+K)}{|s|}|x^s|ds\\
&\leq\frac{\left|C_2\log T_m-a\right|^2(A\log^2 T_m+K)}{T_m}x^a
\end{aligned}
\label{eqn514}
\end{equation}
since $|s|\geq T_m$.

Now consider the second integral and the expression
\begin{equation*}
\begin{aligned}
\log\zeta(s)=\int_{s_c}^{s}\frac{\zeta'(s)}{\zeta(s)}+\log \zeta(s_c).
\end{aligned}
\end{equation*}
Therefore the bound on $\log\zeta(s)$ over this segment becomes

\begin{equation*}
\begin{aligned}
\Bigg.|\log\zeta(s)|\Bigg|_{(s_m,s_c)}&\leq\int_{s_c}^{s} A\log^2 T_m ds+|C_2 \log T_m-a|(A \log^2 T_m+K) \\
&<|s-s_c|\Big\{(A \log^2 T_m+|C_2 \log^2 T_m-a|(A \log^2 T_m+K)\Big\}
\end{aligned}
\end{equation*}
by  \eqref{eqnlog1}.  Now, for $s\in (s_m,s_c)$ we have $|s-s_c|<T_m$, then

\begin{equation}
\Bigg.|\log\zeta(s)|\Bigg|_{(s_m,s_c)}\leq T_m \Big\{A\log^2 T_m+|C_2 \log^2 T_m-a|(A \log^2 T_m+K)\Big\},
\label{logmax}
\end{equation}
and therefore

\begin{equation}
\begin{aligned}
\left|\int_{s_c}^s \frac{\log \zeta(s)}{s} x^s\right|&\leq |s-s_c|\frac{T_m\Big\{ A\log^2 T_m+|C_2 \log^2 T_m-a|(A\log^2 T_m+K\Big\}}{|s|}|x^s|ds \\
&\leq \frac{T_m^2\Big\{ A\log^2 T_m+|C_2\log^2 T_m-a|(A\log^2 T_m+K)\Big\}}{T_m}\frac{1}{T_m^3}.
\end{aligned}
\label{eqn516}
\end{equation}
This is because for $\sigma<\sigma_c$, $|x^s|<(T_m)^{-3}$.  Combining \eqref{eqn514} and \eqref{eqn516} gives
\begin{equation*}
\begin{aligned}
\left|\hspace{2pt}\mathop\int\limits_{\text{Gray($T_m$)}}
\frac{\log\zeta(s)}{s}x^s ds\right|
&\leq\frac{|C_2\log T_m-a|^2(A \log^2 T_m+K)}{T_m} x^a \\
&\hspace{20pt}+\frac{A\log^2 T_m+|C_2 \log^2 T_m-a|(A \log^2 T_m+K)}{T_m^2},
\end{aligned}
\end{equation*}
and as $T_m\to\infty$, the quantity on the left side tends to zero.  The same argument then applies to the lower contour at $-T_m$.
%
% Gray contour over vertical legs
%
\subsubsection{Integral over vertical gray legs}
We wish to first place upper bounds on the quantity $\log\zeta(s)$ over the contours $v_n$ shown at the left in Figure (\ref{fig:file1}).  For $v_1$ we can write 

\begin{equation*}
\Bigg.\log\zeta(s)\Bigg |_{v_1}=\int_{-T_m+i T_m}^{s}\frac{\zeta'(s)}{\zeta(s)}ds+\log\zeta(-T_m+iT_m)
\end{equation*}
with $s\in v_1$.  Using  \eqref{logmax} then:

\begin{equation}
\begin{aligned}
\Bigg.|\log\zeta(s)|\Bigg|_{v_1}&\leq\int_{-T_m+iT_m}^{s}\frac{\zeta'(s)}{\zeta(s)}ds+T_m \Big\{A\log^2 T_m+|C_2 \log^2 T_m-a|(A \log^2 T_m+K)\Big\}\\
&\hspace{-25pt}\leq |-T_m+iT_m-s_1|A\log(|s|+1)+T_m \Big\{A\log^2 T_m+|C_2 \log^2 T_m-a|(A \log^2 T_m+K)\Big\}\\
&\hspace{-25pt}\leq |-T_m+iT_m-s_1|A\log(2T_m+1)+T_m \Big\{A\log^2 T_m+|C_2 \log^2 T_m-a|(A \log^2 T_m+K)\Big\}\\
&\hspace{-25pt}\leq T_m A\log(2T_m+1)+T_m\Big\{A\log^2 T_m+|C_2 \log^2 T_m-a|(A \log^2 T_m+K)\Big\}\\
&\hspace{-25pt}\leq T_m\Big\{A\log(2T_m+1)+A\log^2 T_m+|C_2 \log^2 T_m-a|(A \log^2 T_m+K)\Big\}.
\end{aligned}
\label{eqn900}
\end{equation} 

For the contour $v_2$ over the vertical interval $(s_1, s_2)$ we have
\begin{equation}
\Bigg.\log\zeta(s)\Bigg |_{v_2}=\int_{s_1}^{s}\frac{\zeta'(s)}{\zeta(s)}ds+\log\zeta(s_1).
\end{equation}
Now, $v_1$ is separated from $v_2$ by a branch cut with
\begin{equation*}
\Bigg.\log\zeta(s_1)\Bigg|_{v_2}=\Bigg.\log\zeta(s_1)\Bigg|_{v_1}-2\pi i.
\end{equation*}
We can then write 
\begin{equation*}
\begin{aligned}
\Bigg.\log\zeta(s)\Bigg|_{v_2}=\int_{s_1}^{s}\frac{\zeta'(s)}{\zeta(s)}ds+\Bigg.\log\zeta( s_1)\Bigg|_{v_1}-2\pi i,
\end{aligned}
\end{equation*} 
and using  \eqref{eqn900} we obtain

\begin{equation*}
\begin{aligned}
\Bigg.|\log\zeta(s)|\Bigg|_{v_2}& \\
&\hspace{-40pt}\leq T_m A\log(2T_m+1)+T_m\Big\{A\log(2T_m+1)+A\log^2 T_m+|C_2 \log^2 T_m-a|(A \log^2 T_m+K)\Big\}-2\pi \\
&\hspace{-40pt}\leq T_m\Big\{2A\log(2T_m+1)+A\log^2 T_m+|C_2 \log^2 T_m-a|(A \log^2 T_m+K)\Big\}-2\pi,
\end{aligned}
\end{equation*}
and in general:
\begin{equation}
\begin{aligned}
\Bigg.|\log\zeta(s)|\Bigg|_{v_n}& \\
&\hspace{-50pt}\leq T_m A\log(2T_m+1)\\
&\hspace{-40pt}+T_m\Big\{(n-1)A\log(2T_m+1)+A\log^2 T_m+|C_2 \log^2 T_m-a|(A \log^2 T_m+K)\Big\}+2n\pi \\
&\hspace{-50pt}\leq T_m\Big\{nA\log(2T_m+1)+A\log^2 T_m+|C_2 \log^2 T_m-a|(A \log^2 T_m+K)\Big\}+2n\pi. \\
\end{aligned}
\label{eqn658}
\end{equation}

Now denote the number of non-trivial zeros of $\zeta(s)$ in the range $0<\text{Im}(s)<T$ by $N(T)$.  A proof of the following theorem, first proposed by Riemann, can be found in Ingham \cite{Ingham}:
\begin{theorem}
When $T\to\infty$,
\begin{equation*}
\begin{aligned}
N(T)&=\frac{T}{2\pi}\log\frac{T}{2\pi}-\frac{T}{2\pi}+A\log T.\\
\end{aligned}
\end{equation*}
\label{zerolimit}
\end{theorem}
Thus in the limit as $T_m\to\infty$, the number of zeros in the range $0<t<T_m$ and thus the number of contours $v_n$ does not have an order which exceeds some constant times the factor of $T_m\log T_m$.  Writing  \eqref{eqn658} as
\begin{equation*}
\begin{aligned}
\Bigg.|\log\zeta(s)|\Bigg|_{v_n}\\
&\hspace{-30pt}\leq T_m\Big\{N(T_m)A\log(2T_m+1)+A\log^2 T_m+|C_2 \log^2 T_m-a|(A \log^2 T_m+K)\Big\}\\
&\hspace{-10pt}+2N(T_m)\pi,
\end{aligned}
\end{equation*} 
we have

\begin{equation*}
\begin{aligned}
\left|\int_{s_n}^{s_{n+1}}\frac{\log\zeta(s)}{s}x^s ds\right|\\
&\hspace{-30pt}\leq \frac{|s_n-s_{n+1}| T_m\Big\{N(T_m)A \log(2T_m+1)+G(\log^2 T_m)\Big\}+2N(T_m)\pi}{|s|}|x^s|.
\end{aligned}
\end{equation*}
Since:
\begin{equation*}
\begin{aligned}
N(T)&=\frac{T}{2\pi}\log\frac{T}{2\pi}-\frac{T}{2\pi}+A\log T\\
&\leq\text{O}\left(T\log T\right)\\
&\leq C_3 T\log(T),
\end{aligned}
\end{equation*}
and
\begin{equation*}
\begin{aligned}
|s_n-s_{n+1}|&<T_m \\
|s|&>T_m,
\end{aligned}
\end{equation*}
this gives:

\begin{equation*}
\begin{aligned}
\lim_{T_m\to\infty}
\left|\int_{s_n}^{s_{n+1}}\frac{\log\zeta(s)}{s}x^s ds\right|&\\
&\hspace{-100pt}\leq\lim_{T_m\to\infty} \frac{T_m^2\Big\{T_m\log T_mA \log(2T_m+1)+G(\log^2 T_m)\Big\}+C_3 T_m \log(T_m) \pi}{T_m}\left(\frac{1}{T_m}\right)^3.
\end{aligned}
\end{equation*}
And as $T_m\to\infty$, this limit approaches zero.  Summing over all the contours $v_n$ as the boundary extends to infinity we obtain:

\begin{equation*}
\lim_{m\to \infty}\sum_{n=1}^{N(T_m)+1}\int_{s_n}^{s_{n+1}}\frac{\log\zeta(s)}{s}x^s ds\to 0
\end{equation*}

A similar argument then applies to the lower set of vertical contours.

%
% section
%
\section{Assembling the parts}
We now have all the contours of B(m) and can write  \eqref{eqn002} as  

\begin{equation}
f(x)=P.V. \int _{-2}^1\frac{x^s}{s} ds-\sum_{\textnormal{Im}(\rho_n)>0}\Big[\textnormal{Li}(x^{\rho_n})+\textnormal{Li}(x^{1-\rho_n})\Big]+\sum_{n=2}^{\infty }  \Gamma [0,2n \log (x)]-\text{Log} (2).
\label{eqn62}
\end{equation}
Riemann obtained

\begin{equation}
f(x)=\textnormal{Li}(x)-\sum _{\textnormal{Im}(\rho )>0} \left[\text{Li}\left(x^{\rho }\right)+\text{Li}\left(x^{1-\rho }\right)\right]+\int_x^{\infty } \frac{1}{t\left(t^2-1\right) \text{Log}(t)} \, dt-\text{Log}(2).
\label{eqn003b}
\end{equation}
Therefore, if we can show:
$$
\text{Li}(x)+\int_x^{\infty } \frac{1}{t\left(t^2-1\right) \Log (t)} \, dt=P.V. \int _{-2}^1\frac{x^s}{s}\text{ds}+\sum _{n=2}^{\infty }  \Gamma [0,2n \Log (x)],
$$
then the results calculated via the Residue Theorem are equal to that obtained by Riemann.

Note first if we let $u=x^s$ in the integral
$$
P.V. \int _{-2}^1\frac{x^s}{s}\text{ds},
$$
we obtain
$$
P.V. \int _{-2}^1\frac{x^s}{s}\text{ds}=\text{Li}(x)-\int_0^{1\left/x^2\right.} \frac{1}{\Log  u} \, du.
$$
This leaves us to verify:
$$
\int_x^{\infty } \frac{1}{t\left(t^2-1\right) \Log (t)} \, dt=\sum _{n=2}^{\infty }  \Gamma [0,2n \Log (x)]-\int_0^{1\left/x^2\right.} \frac{1}{\Log  u} \, du.
$$
Now
$$
\Gamma (0,2n \Log x)=\int _{2n \Log  x}^{\infty }t^{-1} e^{-t}dt,
$$
and
$$
\int_x^{\infty } \frac{1}{t\left(t^2-1\right) \Log (t)} \, dt=\sum _{n=1}^{\infty } \int_x^{\infty } \frac{t^{-2n-1}}{\Log  t} \, dt.
$$
We therefore wish to show
$$
\sum _{n=1}^{\infty } \int_x^{\infty } \frac{t^{-2n-1}}{\Log  t} \, dt=\sum _{n=2}^{\infty } \int _{2n \Log  x}^{\infty }t^{-1} e^{-t}dt-\int_0^{1\left/x^2\right.} \frac{1}{\Log  u} \, du.
$$
We can write
$$
\int_x^{\infty } \frac{t^{-2n-1}}{\Log  t} \, dt=\int _{2n \Log  x}^{\infty }t^{-1}e^{-t}dt.
$$
We are left then with showing
$$
\int _{2 \Log  x}^{\infty }t^{-1}e^{-t}dt=-\int_0^{1\left/x^2\right.} \frac{1}{\Log  u} \, du,
$$
which we can show by making two change of variables. The first let $v=\Log u$ to obtain
$$
\int_0^{1\left/x^2\right.} \frac{1}{\Log  u} \, du=\int_{-\infty }^{-2\Log  x} \frac{e^v}{v} \, dv.
$$
Now let v=-r:
$$
\int_{-\infty }^{-2\Log  x} \frac{e^v}{v} \, dv=-\int _{2 \Log  x}^{\infty }r^{-1}e^{-r}dr
$$
which was to be shown. 

\section{Conclusions}

Using a particular holomorphic branch of $\log\zeta(s)$ and the Residue Theorem, we obtain the same expression for $f(x)$ as Riemann.

\vspace{20pt}
\address{\noindent D.C. Milioto\\1178 Hwy. 678\\ LaPlace, LA. 70068 (U.S.A) \\1-985-651-9726 \\ email: miliotodc@rtconline.com}

\end{document}